\documentclass{amsart}

\usepackage{amssymb}
\usepackage[all]{xy}

\newtheorem{thm}{Theorem}

\newtheorem{lem}[thm]{Lemma}
\newtheorem{cor}[thm]{Corollary}

\newtheorem{prop}[thm]{Proposition}

\theoremstyle{definition}
\newtheorem{defn}[thm]{Definition}

\newtheorem{say}[thm]{}
\newtheorem{exmp}[thm]{Example}

\newtheorem*{ack}{Acknowledgments}      
\newtheorem{notation}[thm]{Notation}   
  
\newtheorem{defn-thm}[thm]{Definition--Theorem}  
\newtheorem{defn-lem}[thm]{Definition--Lemma}  

\newtheorem{comments}[thm]{Comments}

\theoremstyle{remark}

\let \cedilla =\c
\renewcommand{\c}[0]{{\mathbb C}}  

\renewcommand{\o}[0]{{\mathcal O}} 
\newcommand{\z}[0]{{\mathbb Z}}

\renewcommand{\a}[0]{{\mathbb A}}

\newcommand{\p}[0]{{\mathbb P}}
\newcommand{\f}[0]{{\mathbb F}}
\newcommand{\q}[0]{{\mathbb Q}}
\newcommand{\map}[0]{\dasharrow}
\newcommand{\qtq}[1]{\quad\mbox{#1}\quad}
\newcommand{\spec}[0]{\operatorname{Spec}}
\newcommand{\pic}[0]{\operatorname{Pic}}

\newcommand{\pico}[0]{\operatorname{\mathbf{Pic}}^{\circ}}

\newcommand{\alb}[0]{\operatorname{Alb}}

\newcommand{\ns}[0]{\operatorname{NS}}

\newcommand{\gal}[0]{\operatorname{Gal}}

\newcommand{\supp}[0]{\operatorname{Supp}}    
\newcommand{\red}[0]{\operatorname{red}}    
\newcommand{\codim}[0]{\operatorname{codim}}    
\newcommand{\im}[0]{\operatorname{im}}

\newcommand{\coker}[0]{\operatorname{coker}}    
    
\newcommand{\Hom}[0]{\operatorname{Hom}}

\newcommand{\sing}[0]{\operatorname{Sing}}

\newcommand{\chow}[0]{\operatorname{Chow}}
\newcommand{\chr}[0]{\operatorname{char}}

\newcommand{\cl}[0]{\operatorname{Cl}}
\newcommand{\clns}[0]{\operatorname{\cl^{\rm ns}}}
\newcommand{\clo}[0]{\operatorname{\cl^{\circ}}}

\newcommand{\rdown}[1]{\lfloor{#1}\rfloor}

\newcommand{\onto}[0]{\twoheadrightarrow}

\newcommand{\lcm}[0]{\operatorname{lcm}}

\newcommand{\isom}[0]{\operatorname{Isom}}

\newcommand{\depth}[0]{\operatorname{depth}} 
\newcommand{\tsum}[0]{\textstyle{\sum}}

\newcommand{\wdiv}[0]{\operatorname{WDiv}} 
 
\newcommand{\pdiv}[0]{\operatorname{PDiv}}

\def\into{\DOTSB\lhook\joinrel\to}

\def\loccoh#1.#2.#3.#4.{H^{#1}_{#2}(#3,#4)}

\DeclareMathAlphabet{\mathchanc}{OT1}{pzc}%
                                {m}{it}

\usepackage[all]{xy}\xyoption{dvips}

\newcommand{\norm}[0]{\operatorname{norm}}

\newcommand{\mdiv}[0]{\operatorname{MDiv}}
\newcommand{\mcl}[0]{\operatorname{MCl}}

\newcommand{\dsh}[0]{\operatorname{DSh}}

\newcommand{\sima}[0]{\stackrel{\rm alg}{\sim}}

\begin{document}
\bibliographystyle{amsalpha}

 \title{Mumford divisors}
 \author{J\'anos Koll\'ar}

\begin{abstract} We define the notion of  Mumford divisors, argue that they are  the natural divisors to study on reduced but non-normal varieties and prove a structure theorem for the Mumford class group. 
\end{abstract}

 \maketitle

In higher dimensional algebraic geometry it is necessary to work with pairs $(X, \Delta)$ where $X$ is not normal. This happens especially  frequently in moduli problems and in proofs that use  induction on the dimension.
While the main interest is in semi-log-canonical  pairs,
for a general theory the right framework seems to be  demi-normal 
 and, more generally, seminormal varieties; see \cite[Sec.5.1]{kk-singbook}.
The aim of this note is to discuss the basics of the divisor theory on
such non-normal schemes.
The main result, Theorem~\ref{main.thm}, shows that the divisor class group of 
  seminormal varieties is much larger than for  normal varieties.

For singular curves a theory of generalized Jacobians was worked out by Severi \cite{MR0024985} and Rosenlicht
 \cite{MR0061422}; an exposition is given by Serre \cite{MR0103191}. 
This theory starts with a smooth, projective  curve $C$ and   a finite subcheme $P\subset C$. The basic objects are 
divisors $B$ that are contained in the open curve  $ C\setminus P$, and two divisors
$B', B''$ are called   linearly equivalent modulo $P$ if  there is a rational function $\phi$ on $C$  such that $(\phi)=B'-B''$ and $1-\phi$ vanishes on $P$.
The linear equivalence classes of degree 0 divisors form a semi-abelian variety. 

This approach does not generalize well to higher dimensions, since the condition that $\phi|_P$ be constant is too stringent if $\dim P>0$. 

A more direct precursor of our definition  is  Mumford's observation that in order to get a good theory of pointed curves  $(C, P)$,
the curve should be nodal and  the marked points $P=\{p_1,\dots, p_n\}$ should  be smooth points of $C$.  Two such point sets $P', P''$ are then considered linearly equivalent iff the corresponding sheaves  $\o_C(P')$ and $\o_C(P'')$ are isomorphic. Note that these sheaves are locally free since the marked points are smooth. 

For moduli purposes, the best  higher dimensional generalization of pointed curves is the class of semi-log-canonical pairs $(X, \Delta)$.
Here $X$ is allowed to have normal crossing singularities in codimension 1, but none of the irreducible components of $\supp\Delta$ is contained in $\sing X$. Equivalently, $X$ is smooth at the generic points of $\supp\Delta$. The latter is a key feature of such pairs and leads to our concept  of Mumford divisors in Definition~\ref{mumf.div.defn}.

Another important property of semi-log-canonical pairs is that
$K_X+\Delta$ is $\q$-Cartier. However $\Delta$ need not be $\q$-Cartier and
even if both $K_X+\Delta$ and $\Delta$ are $\q$-Cartier, the
log-pluricanonical sheaves $\omega_X^{[m]}(\rdown{m\Delta})$ are frequently not.
 For me the immediate reason for studying Mumford divisors  was to understand  
the deformation invariance of log-plurigenera of 
 semi-log-canonical pairs  \cite{k-lpg1, k-lpg2}.
 
If $X$ is normal then there are ways to reduce the general setting to the $\q$-Cartier case by a small modification  \cite[Cor.21]{k-lpg1}, but in the non-normal case this is not possible. Thus it seems necessary to study the divisor theory of semi-log-canonical varieties without any $\q$-Cartier conditions.
This has also long been the viewpoint advocated by Shokurov; see
\cite{sho-3ff,   ambro, fuj-book}.

\section{Mumford divisors}

Recall that a  {\it   Weil divisor}   on a scheme $X$ is  a  formal linear combination  $B=\sum_i m_i B_i$  of codimension 1, irreducible, reduced subschemes   $B_i\subset X$.  The set of all
Weil divisors with coefficients in a ring $R$, also called  {\it   Weil $R$-divisors}, forms  an $R$-module, denoted by   $\wdiv_R(X)$. 

We also write $\wdiv(X)$ to denote $\wdiv_{\z}(X)$.

Combining the notion of Weil divisors with Mumford's observation
gives the following   definition in  higher dimensions.

\begin{defn} \label{mumf.div.defn}
Let $X$ be a scheme.  
A {\it Mumford divisor} on $X$ is a  Weil divisor $B$  such that $X$ is regular at all generic points of $\supp B$. The set of all
Mumford divisors with coefficients in a ring $R$, also called  {\it   Mumford $R$-divisors},  forms an $R$-module, denoted by   $\mdiv_R(X)$. We also write $\mdiv(X)$ to denote $\mdiv_{\z}(X)$.
It is a subgroup of $\wdiv(X)$.

On a normal scheme every  Weil divisor is a Mumford divisor, thus the latter notion is of interest mainly on non-normal schemes. 

More generally, let $X$ be a scheme and $D=\{D_1,\dots, \}$  a set of codimension 1 irreducible, reduced  subschemes  $D_j\subset X$ 
such that $X$ is regular at all codimension 1 points of $X\setminus D$.
A {\it   Mumford divisor}   on $(X, D)$ is  a  Weil divisor $B=\sum_i m_i B_i$ on $X$  such that none of the $B_i$ is among the $D_j$.  As before, the set of all
Mumford divisors with coefficients in a ring $R$ forms an $R$-module, denoted by   $\mdiv_R(X,D)$. It can be viewed as a subgroup of $\wdiv_R(X\setminus D)$, though this is not very useful. 

Thus a Mumford divisor on $X$ is the same as a 
  Mumford divisor on the pair $(X, \sing X) $.
(This is somewhat sloppy notation. We should consider
 the pair $\bigl(X, D(\sing X)\bigr)$ where $D(\sing X)$ is the set of closures of all codimension 1 singular points of $X$.)
Thus $\mdiv(X)=\mdiv(X, \sing X)$.

While the definition makes sense if $X$ is arbitrary and even if $D$ is infinite, in order to avoid various pathologies 
 we only consider schemes  $X$  that satisfy the following.

\medskip
{\it Assumption \ref{mumf.div.defn}.1.}  $X$ is  noetherian, reduced, separated, pure dimensional  and $\sing X$ is a nowhere dense, closed subscheme.
(Note that if $X$ is also excellent then 
it satisfies  the latter  assumption on $\sing X$.)

\medskip

Two    Mumford divisors $B', B''$   on $(X, D)$ are called {\it linearly equivalent modulo $D$} if there is a rational function $\phi$ such that  $(\phi)=B'-B''$ and $\phi$ is regular (and nonzero) at the generic points of $\eta_j\in D$.
 The linear equivalence classes of  Mumford divisors form a group, denoted by $\mcl(X,D)$. It is called the
{\it Mumford class group} of $(X, D)$. The  Mumford class group of $(X, \sing X)$ is called the {\it Mumford class group} of $X$ and denoted by
$\mcl(X)$.  

\medskip
{\it Warning \ref{mumf.div.defn}.2.} 
In the usual definition of  rational equivalence of divisors and  $\chow_{n-1}(X)$, one allows functions that are not regular at the generic points of $\sing X$, cf.\  \cite[Sec.1.2]{Fulton84}.  This implies that rational equivalence is preserved by push-forward  \cite[Sec.1.4]{Fulton84}.
Therefore, if $X$ is a proper variety of dimension $n$ with normalization $\pi:\bar X\to X$, then  $\pi_*:\chow_{n-1}(\bar X)\to \chow_{n-1}(X)$ has
finitely generated kernel and finite cokernel.

By contrast, linear equivalence modulo $D$ is {\em not} preserved by push-forward and
we will see that the difference between  $\mcl(X)$ and
 $\mcl(\bar X)\cong \cl(\bar X)$ is huge. 

Thus the  condition that $\phi$ be regular  at the generic points of $\eta_j\in D$ turns out to have a decisive influence on the structure of $\mcl(X)$. 
\end{defn}

Our main interest is to study $\mcl(X)$  for demi-normal schemes.
The main result, Theorem~\ref{main.thm}, shows that if $\dim X\geq 2$ and $X$ is not normal, then $\mcl(X)$ is much larger than class groups of normal varieties. In particular, if we work over $\c$, then
the discrete part of   $\mcl(X)$ always has uncountably
infinite rank. First we note that
on a normal variety, we do not get anything new.

\begin{lem}\label{norm.mum=weil.lem} Let  $X$ be a   noetherian, normal scheme and $D$ a finite set of divisors. Then $\mcl(X, D)\cong \cl(X)$.
\end{lem}

Proof.  By Lemma~\ref{Tag.09NN.1} every  Weil divisor $B$ is linearly equivalent to  Weil divisor $B'$ whose support does not contain any of the $D_j$. Thus the natural map 
$  \mcl(X, D)\to \cl(X)$ is surjective and it is clearly an injection. \qed
\medskip

Extending the correspondence between Cartier divisors and line bundles to
Mumford divisors leads to divisorial sheaves.

\begin{defn}[Divisorial sheaves] A  coherent sheaf $F$  an a scheme $X$ is called a {\it divisorial sheaf} if there is a closed subset $Z\subset X$ (depending on $F$) such that
\begin{enumerate}
\item $\depth_ZF\geq 2$,
\item $F|_{X\setminus Z}$ is locally free of rank 1 and 
\item  $X\setminus Z$ is $S_2$.
\end{enumerate}
Let $j:U\into X$ be a dense, open, $S_2$  subscheme and $L_U$ an  invertible sheaf on $U$. If $j_*(L_U)$ is coherent then it is a divisorial sheaf. 

If $X$ is pure dimensional and excellent, then $j_*(L_U)$ is coherent iff
$\codim_X(X\setminus U)\geq 2$. This is the only case that we use; see
\cite{k-coherent} for the general setting.

The divisorial sheaves form a group, where the group operation is the
$S_2$-hull of the tensor product.  This is obvious from the second definition in  the  pure dimensional,  excellent case. In general it follows from
\cite{k-coherent}. This group is denoted by $\dsh(X)$.

\end{defn}

\begin{prop} Let $X$ be a  noetherian, reduced, separated, pure dimensional, excellent scheme.  Then $B\mapsto \o_X(B)$ gives an isomorphism
$$
\mcl(X)\cong \dsh(X).
$$
\end{prop}

Proof. Let  $B$ be  a Mumford divisor on $X$.  Then $Z:=\supp B\cap \sing X$
is a closed subset of codimension $\geq 2$ and $B|_{X\setminus Z}$ is Cartier.
Thus $$\o_X(B):=j_*\bigl(\o_{X\setminus Z}(B|_{X\setminus Z})\bigr)$$ is a divisorial sheaf, where $j: {X\setminus Z}\into X$ is the natural open  embedding.
This defines $\mcl(X)\to \dsh(X)$  and it is clearly an injection.

Let $L$ be a divisorial  sheaf and $x_1, \dots, x_r$ the generic points of $\sing X$ that have codimension 1 in $X$. By Lemma~\ref{Tag.09NN.1}  $L$ has a rational  section $s$ that is regular and nonzero  at the points $x_1, \dots, x_r$. Then $(s)$ is a Mumford divisor and
$\o_X((s))\cong L$.  Thus the map is also surjective. \qed

\begin{lem} \label{Tag.09NN.1}
Let $X$ be a reduced, separated, noetherian scheme, $z_1,\dots,z_r$ codimension 1 points of $X$, $x_1,\dots,x_s$ codimension 1 regular points of $X$ and
$m_1,\dots, m_s\in \z$. Let  
$L$ be a line bundle on $X$. Then $L$ has a rational section $\phi$ that is regular at the points  $z_1,\dots,z_r$ and vanishes to order $m_1,\dots, m_s$ at the points $x_1,\dots,x_s$.
\end{lem}

Proof. By  \cite[Tag 09NN]{stacks-project} there is an open affine subscheme
$U\subset X$ that contains  the points  $z_1,\dots,z_r$ and $x_1, \dots, x_s$.
We can thus assume that $X$ is affine.  We can then choose a global section $s\in H^0(X, L)$ that does not vanish at any of  the points  $z_1,\dots,z_r$ and $x_1, \dots, x_s$ and regular functions $\phi_1,\dots,\phi_s$ such that
$\phi_i$ vanishes at $x_i$ to order 1 and does not vanish at any of the other points. Then  $\phi:=s\cdot \prod_i \phi_i^{m_i}$ works. \qed

\begin{say}[$S_2$-hulls]\label{S2.hull.say}  Let $X$ be a reduced, noetherian, affine scheme such that $\sing X$ is closed. Let $g=0$ be an equation of $\sing X$ that is nonzero at all generic points. Let $p_i\in X$ be the embedded primes of $\o_X/(g)$ and $Z\subset X$ the union of their closure.  Then
$\codim_XZ\geq 2$ and $X\setminus Z$ is $S_2$. We call $Z$ the
{\it non-$S_2$-locus} of $X$ and $X\setminus Z$ the 
{\it $S_2$-locus} of $X$. 

Let $j:X\setminus Z\into X$ be the natural embedding. If
$j_*\o_{X\setminus Z}$ is coherent (for example, $X$ is excellent)
then set $X^*:=\spec_Xj_*\o_{X\setminus Z}$. 
Thus $X^*$ is $S_2$, the projection $\pi: X^*\to X$ is finite and an isomorphism outside subsets of codimension $\geq 2$. This $X^*$ is unique and called the
{\it  $S_2$-hull} of $X$. 

For a reduced, noetherian, excellent scheme $X$ we can glue the local
 $S_2$-hulls together to get the  $S_2$-hull  $\pi:X^*\to X$.
If $X$ is seminormal then so is $X^*$.

Since the notion of Mumford divisors is not sensitive to
closed subsets of codimension $\geq 2$ we see that $\pi_*$ induces isomorphisms
$$
\mdiv(X^*)\cong \mdiv(X)\qtq{and} \mcl(X^*)\cong \mcl(X).
\eqno{(\ref{S2.hull.say}.1)}
$$
\end{say}

\begin{say}[Pulling back  divisors]\label{mumf.div.pull.b}
Let $g:Y\to X$ be a morphism and $B$ a Weil divisor on $X$ that is Cartier on an open set $U\subset X$. Then $g^*\bigl(B|_U\bigr)$ is a Cartier divisor on $g^{-1}(U)$. If $Y\setminus g^{-1}(U)$ has  codimension $\geq 2$ then
the closure of $g^*\bigl(B|_U\bigr)$ is a well defined Weil divisor on $Y$, denoted by $g^*(B)$. (Probably $g^{[*]}(B)$ would be a better notation.)

Thus we need to understand morphisms $g:Y\to X$  for which
$\codim_XZ\geq 2$ $\Rightarrow$ $\codim_Yg^{-1}(Z)\geq 2$. 
To guarantee this, we consider the following setting.

\medskip

{\it Assumption \ref{mumf.div.pull.b}.1.}  $X$ is pure dimensional,  the morphism  $g$ is finite and  {\it dominant on  irreducible components.} 
That is, every irreducible component of $Y$ dominates some irreducible component of $X$.
\medskip

 If, in addition,  $X$ and $Y$  are normal then  the  pull-back is defined on all Weil divisors and we get
$g^*:\wdiv(X)\to \wdiv(Y)$ and $g^*: \cl(X)\to \cl(Y)$.
If $Y$ is normal then the  pull-back is defined on all Mumford divisors and we get
$$
g^*:\mdiv(X)\to \wdiv(Y)\qtq{and} g^*: \mcl(X)\to \cl(Y).
\eqno{(\ref{mumf.div.pull.b}.2)}
$$
In particular, if $\pi:\bar X\to X$ is the normalization of $X$ then we get
$$
\pi^*:\mdiv(X)\to \wdiv(\bar X)\qtq{and} \pi^*: \mcl(X)\onto \cl(\bar X),
\eqno{(\ref{mumf.div.pull.b}.3)}
$$
where the surjectivity follows 
 using Lemma~\ref{Tag.09NN.1}.
\end{say}

\begin{say}[Restriction of Mumford divisors]\label{mumf.div.rest}
Let $S\subset X$  be a closed subscheme.
 $B$ is called a {\it Mumford divisor  along $S$} if 
\begin{enumerate}
\item $\supp B$ does not contain any 
irreducible component of $S$,
\item  $B$ is Cartier  at all generic points of $S\cap \supp B$ and
\item $S$ is regular at all generic points of $S\cap \supp B$.
\end{enumerate}
These imply that $B|_S$ is a well-defined Mumford divisor on $S$.
(It seems that  (1--3) are all necessary for this.)
Furthermore, 
there is a subset $Z\subset S$ of codimension $\geq 2$ such that 
$B$ is Cartier at all points of  $S\setminus Z$. The restriction
 sequence
$$
0\to \o_X(B)(-S)\to \o_X(B)\to \o_S(B|_S)\to 0,
\eqno{(\ref{mumf.div.rest}.4)}
$$
 is left exact everywhere and right exact on $X\setminus Z$.

Assume next that $S\subset X$ is a Cartier divisor. Then
(\ref{mumf.div.rest}.4) gives a natural injection
$$ 
r: \o_X(B)|_S\into \o_S(B|_S),
\eqno{(\ref{mumf.div.rest}.5)}
$$
which is an isomorphism on $S\setminus Z$. Since $\o_S(B|_S) $ is $S_2$ by definition, $r$ is  an isomorphism everywhere
 iff $\o_X(B)|_S$ is $S_2$. This is equivalent to
$\o_X(B)$ being  $S_3$ along $S$. (Recall that a sheaf $F$ is  {\it $S_3$ along}  a closed subscheme $Z$ iff 
  $\depth_xF\geq \min\{3, \dim_xF\}$ holds for every   $x\in Z$.)
 We thus obtain the following observation.
\medskip

{\it Claim \ref{mumf.div.rest}.6.} If $S\subset X$ is a Cartier divisor then the sequence (\ref{mumf.div.rest}.4) is exact iff $\o_X(B)$ is $S_3$ along $S$.\qed
\end{say}

\begin{defn}[Principal divisors]
\label{loc.const.say}
Let $V$ be a normal scheme. The group of {\it principal divisors} on $V$ is denoted by $\pdiv(V)$. It is a subgroup of $\wdiv(V)$. 
The latter group is free by definition, hence so is
$\pdiv(V)$; cf.\ \cite[p.880]{lang-alg}. 
If $V$ is a $k$-scheme of finite type and $\dim V\geq 1$ then the rank of $\pdiv(V)$ is  $|k|$ if $k$ is infinite and countable infinite if $k$ is finite.  We can write this as $|\bar k|$.

Let $V$ be a reduced $k$-scheme.  The  algebraic closure of $k$ in $k(X)$ is the algebra of {\it locally constant functions} on $V$. We denote it by
$k^{\rm lc}(V)$. Thus $k^{\rm lc}(V)\subset k(V)$ is a finite dimensional, reduced $k$-algebra.
Note that both $k(V)$ and $k^{\rm lc}(V)$ are birational invariants.
If $V$ is proper and normal , then
$k^{\rm lc}(V)=k[V]=H^0(V, \o_V)$.

If $V$ is reduced and proper  with normalization $ V^n\to V$  and
$\eta_V\in V$ is the generic point then we have  identifications
$$
\pdiv(V^n)=k(V)^*/k^{\rm lc}(V)^*=k(\eta_V)^*/k^{\rm lc}(\eta_V)^*.
\eqno{(\ref{loc.const.say}.1)}
$$
Next let  $g:W\to V$ be a finite morphism of normal and proper  $k$-schemes that is  dominant on  irreducible components. Let 
$\eta_W\in W$  and $\eta_V\in V$ denote the generic points. Then 
$$
\frac{\pdiv(W^n)}{g^*\pdiv(V^n)}
=\frac{k(W)^*}{\bigl\langle g^*k(V)^*, k^{\rm lc}(W)^*\bigr\rangle}
=\frac{k(\eta_W)^*}{\bigl\langle g^*k(\eta_V)^*, k^{\rm lc}(\eta_W)^*\bigr\rangle}.
\eqno{(\ref{loc.const.say}.2)}
$$
We study the structure of these quotient groups in Paragraph~\ref{mult.ff.say}.
\medskip

{\it Remark \ref{loc.const.say}.3.} It does not seem possible to define a sensible analog of $k^{\rm lc}(V)$ for non-reduced schemes $V$. The problem is that if $\phi$ is a nilpotent global section of $\o_V$ then
$1+\phi$ satisfies the equation  $(x-1)^n=0$ for some $n>0$. Thus there can be too many global sections of  $\o_V$ that satisfy a polynomial
in $k[x]$ with leading coefficent $=1$. 
\end{defn}

\begin{say}[Algebraic equivalence of Mumford divisors]
Let $D\subset X\times C$ be an algebraic family of divisors over a connected curve $C$ such that every fiber $D_c\subset X_c$ is a Mumford divisor. 
By normalization of $X$ we get $D^n\subset \bar X\times C$, which is again
an algebraic family of divisors. Thus 2 Mumford divisors $D_1, D_2$ on $X$ are algebraically equivalent iff their pull-backs  $\pi^*D_1, \pi^*D_2$ are
algebraically equivalent on the normalization $\pi:\bar X\to X$. 
Thus
$$
\mcl(X)/\mcl^{\rm alg}(X)\cong \cl(\bar X)/\cl^{\rm alg}(\bar X)\cong \tsum_i \cl^{\rm ns}(\bar X).
$$
\end{say}

\section{Seminormal and demi-normal schemes}

\begin{say}[Seminormal schemes] \label{demi.norm.say}
A morphism of schemes $p:Y\to X$ is an {\it isomorphism on points} if
$p^*_x: k(x)\mapsto k\bigl(\red p^{-1}(x)\bigr)$ is an isomorphism for every $x\in X$.  If $X,Y$ are $\c$-schemes then this holds iff $p$ is a bijection on $\c$-points.

A  morphism of schemes $p:X'\to X$ is a {\it partial seminormalization}
if $X'$ is reduced, $p$  is finite and an isomorphism on points. A scheme $X$ is called 
 {\it seminormal} if every partial seminormalization $p:X'\to X$ is an isomorphism. 
If the normalization $X^n\to X$ is finite (for example, if $X$ is excellent) then there is  a unique factorization  $X^n\to X^{sn}\to X$ where
$X^{sn}$ is seminormal and  $X^{sn}\to X$ is a partial seminormalization.
Then $X^{sn}\to X$ is called  the 
 {\it seminormalization} of $X$. 

See \cite{MR0239118, MR0266923, MR0277542} for the origin of this notion and \cite[Sec.10.2]{kk-singbook} for more recent details.
\end{say}

\begin{say}[Conductor] \label{conductor.say}
Let $\pi:X'\to X$ be a finite, birational morphism of reduced schemes.  The {\it conductor ideal sheaf} is $I_{X'/X}:=\Hom_X\bigl(\pi_*\o_{X'}, \o_X\bigr)\subset \o_X$. It  is the largest ideal sheaf on $X$ that is also an ideal sheaf on $X'$. 
Somewhat sloppily, I refer to either of the subschemes  $V(I_{X'/X})\subset X$ and $V(I_{X'/X})\subset X'$ as  the conductor of $\pi$. 

The most important special case is when $\pi:\bar X\to X$ is the 
 normalization.
 The conductor ideal  $I_X:=\Hom_X\bigl(\pi_*\o_{\bar X}, \o_X\bigr)\subset \o_X$ of $\pi$ is then  called the {\it conductor ideal} of $X$. 
We write $D:=V(I_X)\subset X$ and $\bar D:=V(I_{\bar X})\subset\bar X$ and call both of them the conductor of $\pi$. 
Note that $\pi_D:\bar D\to D$ is a finite, surjective morphism and  $X$ is uniquely determined by $\bar X$ and the morphism
$\pi_D:\bar D\to D$. Indeed, let $r_{\bar X/\bar D}:\o_{\bar X}\to \o_{\bar D}$ be the restriction map. By push-forward we get
$R_{\bar X/\bar D}:\pi_*\o_{\bar X}\to \pi_*\o_{\bar D}$. 
There is a natural injection $j:\o_D\into \pi_*\o_{\bar D}$ and
$$
\o_X= R_{\bar X/\bar D}^{-1}\bigl(j(\o_D)\bigr),
\eqno{(\ref{conductor.say}.1)}
$$
see \cite[Sec.7.2]{rc-book} for details.

It is easy to see that if $X$ is  seminormal then  $\bar D$ is reduced 
(cf.\ \cite[I.7.2.5]{rc-book}) and if
$X$ is $S_2$ then  $\bar D\subset \bar X$ is pure of codimension 1
(cf.\ \cite[10.14]{kk-singbook}).  
Thus if $X$ is  seminormal and $S_2$ then 
$\pi_D:\bar D\to D$ is a finite, surjective morphism of reduced schemes  that is dominant on  irreducible components. Thus the
 induced map on the function algebras 
$$
\pi_D^*:  k(D)\into k(\bar D),
$$
 is a finite extension of direct sums of fields.
In this case  $X$ is uniquely determined by $\bar X$ and by
$\pi_D^*:  k(D)\into k(\bar D)$, cf.\ \cite[10.14]{kk-singbook}.
\end{say}

\begin{notation}\label{dn.notation}
Let $X$ be a seminormal, $S_2$  scheme with normalization $\pi:\bar X\to X$.  We write $D\subset X$ and $\bar D\subset \bar X$ for its conductors.
We also call $D$ the {\it singular divisor} of $X$. 
The normalization, and its restriction to $\bar D$, are usually written as 
$$
\pi: (\bar X, \bar D)\to (X, D)\qtq{and} \pi_D:=\pi|_{\bar D}: \bar D\to D.
\eqno{(\ref{dn.notation}.1)}
$$
We write 
$$X=\cup_{i\in I}X_i\qtq{and} D=\cup_{j\in J} D_j
\eqno{(\ref{dn.notation}.2)}
$$ as the union of their irreducible components. 
Thus $\bar X=\cup_{i\in I}\bar X_i$ and $\pi_i:\bar X_i\to X_i$ are the normalizations. 

Similarly we have a decomposition  $\bar D=\cup_{j\in J} \bar D_j$ and finite maps  $\pi_j:\bar D_j\to  D_j$. Note that the $\bar D_j$ may be reducible. 

In general, neither the $D_j$ nor the $\bar D_j$ are normal, their
normalizations are denoted by  $D^n_j$ and $\bar D^n_j$. Thus we get finite morphisms
 $\pi_j^n:\bar D_j^n\to  D_j^n$, that are dominant on  irreducible components.
Set  
$$
Z:=\sing D\cup (\mbox{non-$S_2$-locus of $X$})\cup \pi\bigl(\sing \bar X\cup \sing \bar D\bigr)
$$
and define the {\it semi-regular} locus of $(X, D)$ as
$$
(X^{\rm sr}, D^{\rm sr}\bigr):=(X\setminus Z, D\setminus Z).
\eqno{(\ref{dn.notation}.3)}
$$
Note that  $X^{\rm sr}$ is $S_2$, the schemes $\bar X^{\rm sr}, D^{\rm sr}, \bar D^{\rm sr}$ are all regular and $\pi_D^{\rm sr}:\bar D^{\rm sr}\to  D^{\rm sr}$ is flat.
If $\bar X^{\rm sr}, D^{\rm sr}, \bar D^{\rm sr}$ are all smooth, we also call this the {\it semi-smooth} locus. 
\end{notation}

\begin{defn}[Demi-normal schemes]\label{dn.dn.notation} A scheme 
$X$ is called  {\it demi-normal} if it is seminormal, $S_2$ and $\pi_D$ has degree 2 over every
irreducible component of $D$. 
Equivalently, $X$ is $S_2$ and its codimension 1 points are either regular or nodal (see \cite[1.31]{kk-singbook} for nodes in residue characteristic 2).

Using Notation~\ref{dn.notation} we write 
$X=\cup_{i\in I}X_i$ and $D=\cup_{j\in J} D_j$
 as the union of their irreducible components and  $\bar D=\cup_{j\in J} \bar D_j$.  The induced  maps  $\pi_j:\bar D_j\to  D_j$ have degree 2. As before, the $\bar D_j$ may be reducible. 

Let   $\eta_j\in D_j$ and $\bar\eta_j\in\bar  D_j$ denote the generic points.
Thus $k(\bar\eta_j)$ is a reduced $k(\eta_j)$-algebra of dimension 2. 
There are 3 possibilities.
\begin{enumerate}
\item (Reducible)   $k(\bar\eta_j)\cong k(\eta_j)+k(\eta_j)$. Let $\tau_j$ denote the  involution that interchanges the summands. 
\item (Separable)    $k(\bar\eta_j)/k(\eta_j)$ is a degree 2 separable field extension with Galois involution $\tau_j$.
\item (Inseparable)    $k(\bar\eta_j)/k(\eta_j)$ is a degree 2 inseparable field extension. This can happen only if $\chr k(\eta_j)=2$. 
\end{enumerate}
In  cases   (\ref{dn.dn.notation}.1--2)   $\tau_j$ induces an involution on $\bar D_j^n$, which we also denote by $\tau_j$.  All the $\tau_j$ together determine an involution
$\tau: \bar D^n\to \bar D^n$ called the {\it gluing involution.}

Thus $D^n=\bar D^n/(\tau)$ and  $X$ is uniquely determined by the data
$(\bar X, \bar D,\tau)$; see \cite[5.3]{kk-singbook} for details. 
\end{defn}

\section{The structure of the  class group}

Let us  recall the known general results about the class group of normal schemes and of nodal curves.

\begin{defn}[Class group of a normal variety]\label{cl.norm.say}
Let $X$ be a  normal and proper $k$-variety. Let $\clo(X)\subset \cl(X)$ denote the group of divisors that are algebraically equivalent to 0.
I call  the quotient  $\clns(X):=\cl(X)/\clo(X)$ the {\it N\'eron-Severi} class group of $X$. 
Note that there are natural  inclusions
$$
\pic^{\circ}(X)\subset \clo(X)\qtq{and} \ns(X)\subset \clns(X).
\eqno{(\ref{cl.norm.say}.1)}
$$
Both of these are equalities iff every Weil divisor is Cartier, for example if $X$ is regular. 
\end{defn} 

Note that if  $X$ is a 
 proper  $k$-scheme 
with Picard scheme $\pico(X)$ then there is a natural injection
$\pic^{\circ}(X)\into \pico(X)(k)$
which is an isomorphism if $X(k)\neq \emptyset$ and $H^0(X, \o_X)\cong k$, cf.\ \cite[V.2.1]{FGA}.
Furthermore,    $\ns(X)$
is a finitely generated abelian group. In the complex case these go back to  Picard \cite{picard-95} and Severi \cite{MR1511373};
see Mumford's Appendix V in \cite{MR0469915} for a discussion of the (quite convoluted) history of the general case.
The most complete references may be the hard-to-find papers of
Matsusaka   \cite{mats-pic} and  of N\'eron \cite{neron-pic}. The latter  also proves the analogous results for the class groups.

More recent results on various aspects of the class group of singular varieties are discussed in \cite{MR1203911, MR2457299, MR2567426}.

It is quickest to construct the scheme ${\mathbf{Cl}}^{\circ}(X)$ using
the Albanese variety.

\begin{defn}[Albanese variety] Let $X$ be a  geometrically normal and proper $k$-variety. If $X(k)\neq \emptyset$ then, following the classical definition,   the {\it Albanese map} is the universal
rational map  $\operatorname{alb}_X:X\map \alb(X)$ to an Abelian variety. 
 This definition   differs from the one in \cite[VI.3.3]{FGA}, which requires
$\operatorname{alb}_X$ to be a  morphism but assumes universality only for morphisms $X\to (\mbox{Abelian variety})$. 
The 2 versions are the same if $X$ has rational singularities. 

In general one should define
the {\it Albanese map} is the universal
rational map  $\operatorname{alb}_X:X\map \alb(X)$ to an Abelian torsor.  
Here an {\it Abelian torsor} is a principal homogeneous space under an Abelian $k$-variety. 
 (Equivalently, a $k$-scheme $A_k$ such that $A_K$ is isomorphic to an Abelian variety for some separable extension $K/k$. In this case $A_k$ is a principal homogeneous space under $\mathbf{Aut}^{\circ}(A_k)$.)

Over $\c$, their construction of the Albanese variety is  usually attributed to \cite{alb-20, alb-21}
(see \cite{MR1419986}  for a more accessible source), 
but in retrospect much of it is already in \cite{severi-13}. 
Over any field the existence of the Albanese variety is proved by Matsusaka  \cite{mats-pic}.
\end{defn}

We can now state the structure theory of the class group in the following form.

\begin{thm}[Class group of a normal variety]\label{cl.norm.thm}
Let $X$ be a geometrically normal and proper $k$-variety. 
Then 
\begin{enumerate} 
\item there is a natural injection
$\clo(X)\into \pico\bigl(\alb(X)\bigr)(k)$ 
which is an isomorphism if $X(k)\neq \emptyset$, and 
\item   $\clns(X)$   is a finitely generated abelian group. \qed
\end{enumerate}
\end{thm}

If $X$ has a resolution of singularities $Y\to X$ then the following 
can be used to compute $\cl(X)$ in terms of the better known $\pic(Y)$; cf.\ \cite{MR3044492}.

\begin{lem} \label{cl.from.res.lem}
Let $g:Y\to X$ be a proper, birational morphism between normal varieties and $\{E_i: i\in I\}$ the exceptional divisors. Then push-forward by $g$ gives isomorphisms
$$
\clo(Y)\cong \clo(X) \qtq{and}
\clns(Y)/\oplus_i \z[E_i]\cong \clns(X).
\eqno{(\ref{cl.from.res.lem}.1)}
$$
\end{lem}

Proof. The push-forward is clearly surjective, thus we only need to prove the following.
\medskip

{\it Claim \ref{cl.from.res.lem}.2.} If  $D$ is a divisor on $Y$ and $g_*(D)$ is algebraically equivalent to 0 then 
 $D$ is  algebraically equivalent to a sum of exceptional divisors.  
\medskip

If $g_*(D)$ is Cartier then $g^*\bigl(g_*(D)\bigr)-D$ is  exceptional and
algebraically equivalent to $D$. Thus  we claim that
even if $g_*(D)$ is not Cartier, we can still define $g^*\bigl(g_*(D)\bigr)$. 

If $g_*(D)$ is linearly equivalent to 0, then $g_*(D)=(\phi)$ for some rational function  $\phi$ on $X$. Thus $D-(\phi\circ g)$ is $g$-exceptional. 
Consider next a family of  algebraically equivalent
 divisors $D^X\subset X\times C$ and its birational transform
$D^Y\subset Y\times C$ for some irreducible, smooth, projective curve $C$.
Let $0\in C$ be a point and $D_0$ a divisor on $Y$ such that 
 $g_*(D_0)=D^X_0$  and $D_0$ is  algebraically equivalent to a sum of exceptional divisors. Note that $g_*(D^Y_0)=g_*(D_0),$ thus
 $D^Y_0-D_0$ is $g$-exceptional, hence  
  $D^Y_0\sima \tsum a_iE_i$ for some $a_i\in \z$. Therefore, for any other points $c\in C$ we have 
$$
D^Y_c-\tsum a_iE_i\sima D^Y_0-\tsum a_iE_i\sima 0.
$$ 
Such families generate algebraic equivalence, thus (\ref{cl.from.res.lem}.2) holds. 
Finally note  that 
 the classes of exceptional divisors are linearly independent by
\cite[3.39]{km-book}. \qed
\medskip

Combining Lemma~\ref{cl.from.res.lem} and (\ref{cl.norm.thm}.1)
gives the following useful observation.

\begin{cor} Let $X$ be a normal and proper $k$-variety and $k$ perfect.
Then
\begin{enumerate}
\item $\clo(X)$ is a birational invariant of $X$ and
\item there is a normal and proper $k$-variety $X'$ that is birational to $X$ such that $\clo(X')=\pic^{\circ}(X')$. \qed
\end{enumerate}
\end{cor}

\begin{say}[Picard group of nodal curves]\label{over.curve.say}
 Let $C$ be a proper, connected nodal curve over $\c$ and $\pi:\bar C\to C$ its normalization.  Then 
$\mcl(C)\cong \pic(C)$,  $\cl(\bar C)=\pic(\bar C)$ and there is an exact sequence
$$
0\to (\c^*)^m\to \mcl(C)\stackrel{\pi^*}{\longrightarrow} \cl(\bar C)\to 0,
\eqno{(\ref{over.curve.say}.1)}
$$
where $m=h^1(C, \o_C)-h^1(\bar C, \o_{\bar C})$; see \cite[Sec.X.2]{acg2} for details.
More invariantly, 
 there is a topological 1-complex  ${\mathcal D}(C)$ whose points are the irreducible components of $\bar C$ and for each  node $p\in C$  we add an edge that connects  the 2 irreducible components passing through $p$. 
(We add a loop if they are the same.) Then the  improved exact sequence is
$$
0\to  H^1\bigl({\mathcal D}(C), \c^*\bigr)\to \mcl(C)\stackrel{\pi^*}{\longrightarrow} \cl(\bar C)\to 0.
\eqno{(\ref{over.curve.say}.2)}
$$
Our aim is to generalize this sequence to higher dimensions. 
Thus 
let $X$ be a proper, seminormal variety over $\c$
with normalization $\pi: \bar X\to X$. 
Then $\pic(\bar C) $ should be replaced by  $\cl(\bar X)$ and again there is a natural topological 1-complex  ${\mathcal D}_1(X)$, giving
 a subgroup $H^1\bigl({\mathcal D}_1(X), \c^*\bigr) \into \mcl(X)$. 
However,  the sequence
$$
  H^1\bigl({\mathcal D}_1(X), \c^*\bigr)\to \mcl(X)\stackrel{\pi^*}{\longrightarrow} \cl(\bar X)
\eqno{(\ref{over.curve.say}.3)}
$$
is no  longer exact.  
In (\ref{cwtlb.say}--\ref{pwtlb.say}) we define subgroups of the Mumford class group
$$
\begin{array}{lcl}
\mcl^{\rm ct}(X) & := & \ker\bigl[\pi^*:\mcl(X)\to \cl(\bar X)\bigr]\qtq{and}\\
\mcl^{\rm pt}(X) & := & \im\bigl[H^1\bigl({\mathcal D}_1(X), \c^*\bigr)\to \mcl(X)\bigr],
\end{array}
\eqno{(\ref{over.curve.say}.4)}
$$
but the key new problem  is to understand the quotient
$$
\frac{\mcl^{\rm ct}(X)}{\mcl^{\rm pt}(X)}=
\frac{\ker\bigl[\pi^*: \mcl(X)\to \cl(\bar X)\bigr]}{H^1\bigl({\mathcal D}_1(X), \c^*\bigr)},
\eqno{(\ref{over.curve.say}.5)}
$$
which turns out to be essentially  a  free  abelian group of uncountable rank. 
\end{say}

\section{The Mumford class group}

The following result  describes the  Mumford class group  of seminormal varieties over an arbitrary base field. As we noted in Paragraph~\ref{S2.hull.say}, it is enough to consider  $k$-schems that are 
seminormal and  $S_2$.

\begin{thm}\label{main.thm}
 Let $X$ be a proper, geometrically semi-normal, $S_2$,  pure dimensional $k$-scheme with 
normalization $\pi:\bar X\to X$. Following Notation~\ref{dn.notation}, let 
$\{X_i:i\in I\}$ be the irreducible components of $X$,
$\{D_j:j\in J\}$  the irreducible components of the conductor $D\subset X$
and  $\bar D=\cup_{j\in J} \bar D_j$.
Then  the Mumford class group $\mcl(X)$ has a natural filtration by  subgroups
$$
0\subset \mcl^{\rm pt}(X) \subset \mcl^{\rm ct}(X)  \subset \mcl(X),
\eqno{(\ref{main.thm}.1)}
$$
and the successive quotients have the following descriptions.
\begin{enumerate}\setcounter{enumi}{1}
\item  $\mcl^{\rm pt}(X)\cong  \coker\bigl[k^{\rm lc}(\bar X)^*\times k^{\rm lc}(D)^*\rightarrow  k^{\rm lc}(\bar D)^*\bigr]$  as in (\ref{pwtlb.say}.2--3), where $k^{\rm lc} $ denotes  locally constant functions as in 
Definition~\ref{loc.const.say}.
\item $\mcl^{\rm ct}(X)=\mcl^{\rm pt}(X) $ iff either $X$ is normal or $\dim X=1$.
\item If  $\dim X\geq 2$ and  $X$ is not normal  then  $\mcl^{\rm ct}(X)/\mcl^{\rm pt}(X) $  has the following properties.
\begin{enumerate}
\item There is a natural isomorphism
$$
\mcl^{\rm ct}(X)/\mcl^{\rm pt}(X)\cong \oplus_j \  k\bigl( \bar D_j\bigr)^*/
\bigl\langle k\bigl( D_j\bigr)^*, k^{\rm lc}\bigl( \bar D_j\bigr)^*\bigr\rangle.
$$
\item  If $\chr k\neq 2$  then there is a non-canonical isomorphism $$\mcl^{\rm ct}(X)/\mcl^{\rm pt}(X)\cong \z^{|\bar k|}+(\mbox{\rm finite group}).
$$
\item  If $\chr k= 2$  then there is a non-canonical isomorphism 
$$\mcl^{\rm ct}(X)/\mcl^{\rm pt}(X)\cong \z^{|\bar k|}+T_2+(\mbox{\rm finite group}),
$$
where $T_2$ is a $2^m$-torsion group of cardinality $|\bar k|$ for some $m$ and 
 either one (but not both) of the first 2 summands may be missing.
\end{enumerate}
\end{enumerate}
\end{thm}

\begin{comments}\label{main.thm.comms}
A useful feature of these formulas is that each quotient can be computed using
only  small parts of the description of $X$ as the push-out of the maps
$\bar D\to \bar X$ and $\bar D\to D$. 
\begin{enumerate}
\item $ \mcl(X)/\mcl^{\rm ct}(X)$ is computable from 
the normalization $\bar X$,
\item  $ \mcl^{\rm ct}(X)/\mcl^{\rm pt}(X)$ from the generic points of  $\bar D\to D$ and
\item $ \mcl^{\rm pt}(X)$  needs only the generic points of
$\bar X$ and of the conductors.  
\end{enumerate}

Furthermore, we see during the proof that in (\ref{main.thm}.5.b--c) the order of every element of the torsion  summands divides
$\lcm\bigl(\deg (\bar D_j/D_j): j\in J\bigr)$. 

\medskip

{\it Warning \ref{main.thm.comms}.4.}   There does not seem to be any sensible way to identify the quotient $\mcl^{\rm ct}(X)/\mcl^{\rm pt}(X) $ with the $k$-points of an algebraic group.

Also, as far as I can tell,  $\mcl(X)$ does not contain any algebraic subgroup that
maps onto $\oplus_i  \clo\bigl(\bar X_i\bigr)$. 
\end{comments}

\begin{exmp} Set $S:=(xy=0)\subset \p^3$ with double curve  $D=(x=y=0)\cong \p^1$.
 Then $\mcl^{\rm pt}(S)=0$ and $\mcl(S)$  sits an in exact sequence
$$
0\to k(\p^1)^*/k^*\to \mcl(S)\to \cl(\p^2)^2\cong \z^2\to 0.
$$
A free generating set of $k(\p^1)^*/k^* $ is given by the  
$\gal(\bar k/k)$-orbits on  $\p^1(\bar k)$ or as  the set of 
$\gal(\bar k/k)$-orbits on  $\bar k$ plus the point at infinity. 
  Thus 
$$
\mcl(S)\cong \z^{\bar k/\gal(\bar k/k)}+\z^3
$$
is a free abelian group of rank $|\bar k|$.
\end{exmp}

\begin{exmp} Set $S:=(xyz=0)\subset \p^2_{xyz}\times \p^1_{uv}$ with double curves  $D_1, D_2, D_3$ isomorphic to $\p^1_{uv}$.
 Then
\begin{enumerate}
\item  $ \mcl(S)/\mcl^{\rm ct}(S)\cong \cl(\p^1\times \p^1)^3\cong \z^6$.
\item $\mcl^{\rm pt}(S)\cong k^*$.
\item There is an isomorphism
$\mcl^{\rm ct}(S)/\mcl^{\rm pt}(S)\cong \bigl(k(\p^1)^*/k^*\bigr)^3$.
\end{enumerate}
Taken together we get  non-canonical isomorphisms
$$
\mcl(S)\cong k^* + \bigl(k(\p^1)^*/k^*\bigr)^3 + \z^6\cong k^* +\z^{|\bar k|}.
$$
\end{exmp}

\begin{exmp}[Higher pinch point] Set $S:=(x^mw=y^mz)\subset \p^3$ and assume that $\chr k\nmid m$. 
It is singular along the line $D=(x=y=0)$ and has 2 higher pinch points on it at the coordinate vertices. Its normalization is the ruled surface $\f_{m-1}$; the line $(z=w=0)$ gives the negative section. We obtain that
$$
\mcl(S)\cong \z/m\z +\z^{|\bar k|}.
$$
The m-torsion is given by the difference of lines
$B:=(x=z=0)-(y=w=0)$.  Note that
$$
\bigl(\tfrac{z}{w}\bigr)=m(x=z=0)-m(y=w=0)=mB,
$$
so $mB\sim 0$. Furthermore, $\frac{z}{w}=\bigl(\frac{x}{y}\bigr)^m$,
but $\frac{x}{y}$ is not regular along the singular line $D$, so
it does not give a linear equivalence for $B$.
\end{exmp}

 Next we define the two subgroups that appear  in (\ref{main.thm}.1). 

\begin{say}[Component-wise trivial Mumford divisors]\label{cwtlb.say}
Let $X$ be a reduced scheme with normalization
$\pi:\bar X\to X$ and conductor $\bar D\subset \bar X$. 
As in (\ref{mumf.div.pull.b}.3) we have a natural pull-back map on divisors 
$$
\pi^*: \mdiv(X)\onto \mdiv(\bar X, \bar D),
\eqno{(\ref{cwtlb.say}.1)}
$$
which descends to a pull-back map on divisor classes
$$
\pi^*: \mcl(X)\onto \mcl(\bar X, \bar D)\cong \cl(\bar X),
\eqno{(\ref{cwtlb.say}.2)}
$$
where the last isomorphism is by Lemma~\ref{norm.mum=weil.lem}.
The kernel of this map is called the group of
{\it component-wise trivial} Mumford divisors. We denote it by
$\mcl^{\rm ct}(X) $. 
Thus there is an exact sequence
$$
0\to \mcl^{\rm ct}(X) \to \mcl(X)\to  \oplus_i  \cl\bigl(\bar X_i\bigr)\to 0
\eqno{(\ref{cwtlb.say}.3)}
$$
and $\mcl(X)/\mcl^{\rm ct}(X)$ is the quotient of $\mcl(X)$ that is easiest to study. 
\end{say}

\begin{say}[Piecewise trivial divisorial sheaves, geometric case]\label{pwtlb.say.1} The geometrically simplest case is when
 $\bar X$ and  $\bar D$ are both smooth and
each $D_j\subset D$ is contained in 2 irreducible components
$X_{i(j,1)}$ and  $ X_{i(j,2)}$. Correspondingly we  get that 
$\bar D_j=\bar D'_j\amalg \bar D''_j$ has 2 irreducible components for every $j$, where $\bar D'_j$ lies on $\bar X_{i(j,1)}$ and  $\bar D''_j$ lies on $\bar X_{i(j,2)}$.
(It is convenient to fix an ordering of the $X_i$ and 
always choosing $i(j,1)<i(j,2)$.)
Thus the gluing involution is given by isomorphisms $\tau_j: \bar D'_j\cong \bar D''_j$. A  {\it piecewise trivial line bundle} on $X$ is given by
the trivial line bundle on $\bar X$ with gluing data
$$
\mu_j: \o_{\bar D'_j}\cong \tau_j^*\o_{\bar D''_j},
\eqno{(\ref{pwtlb.say.1}.1)}
$$
where $\mu_j$  is multiplication by a constant $\mu_j\in k^*$.  
We can change the trivialization of $\o_{\bar X}$ by different multiplicative constants $\nu_i$ on the irreducible components $\bar X_i$. 
This changes $\mu_j$ by the quotient  $\nu_{i(j,1)}/\nu_{i(j,2)}$.

The end result is best expressed in terms of the {\it dual complex} of $X$.
This is a topological 1-complex  ${\mathcal D}_1(X)$ whose points are the irreducible components  $X_i\subset X$ and for each irreducible component   $D_j\subset D$ we add an edge that connects   $X_{i(j,1)}$ and  $ X_{i(j,2)}$.
  (In general we add a loop if an  irreducible component of $X$ has self-intersection along $D_j$.) Then the above considerations say that 
$$
\mcl^{\rm pt}(X)\cong H^1\bigl({\mathcal D}_1(X), k^*\bigr).
\eqno{(\ref{pwtlb.say.1}.2)}
$$

If $X$ is a simple normal crossing scheme or, more generally, a semi-dlt pair, then  ${\mathcal D}_1(X)$  is the 1-skeleton of the true
dual complex of $X$, cf.\ \cite{dkx}.

If  $(X, D)$ is an arbitrary seminormal pair over $\c$ then
in general $\bar X$ and $\bar D$ are not smooth.  However, as we noted in  Notation~\ref{dn.notation}, there is a
closed subset $Z\subset X$ of codimension $\geq 2$ such that 
the smoothness assumptions are satisfies by
$(X\setminus Z, D\setminus Z)$. We can then use the above construction to get
$\mcl^{\rm pt}(X\setminus Z)$ and then pushing  forward from
$X\setminus Z$ to $X$ to obtain 
$$
\mcl^{\rm pt}(X):=\mcl^{\rm pt}(X\setminus Z).
\eqno{(\ref{pwtlb.say.1}.3)}
$$
Note that while the divisors in $\mcl^{\rm pt}(X) $ are Cartier on $X\setminus Z $, they need not be Cartier on $X$. 
\medskip

{\it Example \ref{pwtlb.say.1}.3.} let $S:=(x_1x_2x_3=0)\subset \p^3$ and consider the divisor
$B:=L_1+L_2+L_3-(x_0=0)$ where each $L_i$ is a line through the origin contained in the plane $(x_i=0)$ but different from the coordinate axes. 
Then $B$ is in $\mcl^{\rm pt}(S)$ but it is Cartier only if
the 3 lines $L_1, L_2, L_3$ are coplanar. 

Here we can take $Z=(x_1=x_2=x_3=0)$ and ${\mathcal D}(S\setminus Z)$
is a triangle. Thus 
 $\mcl^{\rm pt}(S)\cong k^*$. The map from the above divisors to
$k^*$  can be given as follows.
Write  $L_i=(x_i=a_ix_{i+1}+b_ix_{i+2}=0)$ using indices modulo 3.  Then 
$$
L_1+L_2+L_3-(x_0=0)\mapsto  \tfrac{a_1a_2a_3}{b_1b_2b_3}\in k^*
$$
gives the isomorphism $\mcl^{\rm pt}(S)\cong k^*$.

\end{say}

\begin{say}[Trivializations]\label{trivialize.say}
Let $X$ be a scheme. 
A {\it trivialization} of a line bundle  $L$ on $X$ is an
isomorphism $\sigma:\o_X\cong L$. 

Let $X$ be a $k$-scheme. 
Two trivializations $\sigma_i:\o_X\cong L$
 are called {\it projectively equivalent} if
$\sigma_2^{-1}\sigma_1:\o_X\to \o_X$ is multiplication by a locally constant function. A   {\it projective trivialization} of a line bundle  $L$
is an equivalence class of projectively equivalent trivializations. 

Note that if $X$ is proper and normal, more generally, if $X$ is the complement of codimension $\geq 2$ subset in a proper and normal $k$-scheme, then any two trivializations are projectively equivalent. 

In our applications this will hold for $X$ but not necessarily for its conductor $D$. The codimension 2 subset $Z\subset X$ that we remove in (\ref{pwtlb.say.1}.3) usually has  codimension 1 in $D$. 
\end{say}

\begin{say}[Piecewise trivial divisorial sheaves, general case]\label{pwtlb.say}
 Let $(X,D)$ be a seminormal $k$-scheme with normalization
$\pi:(\bar X, \bar D)\to (X,D)$. Assume first that $D, \bar D$ are normal and
$\pi_D:\bar D\to D$ is flat.

A {\it piecewise trivial} line bundle on $X$ is a line bundle $L$ on $X$ endowed with  projective trivializations of $L|_D$ and of $\pi^*L$ that are compatible over $\bar D$.

If we fix the trivializations of $\pi^*L$ and of $L|_D$, then 
we are left with locally constant isomorphisms
$$
\isom^{\rm lc}\bigl(\pi^*\o_D, \o_{\bar X}|_{\bar D}\bigr)\cong k^{\rm lc}(\bar D)^*.
\eqno{(\ref{pwtlb.say}.1)}
$$
 Changing the trivializations  of $\pi^*L$ and of $L|_D$  gives
that the group $\pic^{\rm pt}(X)$ of piecewise trivial line bundles on $X$ is isomorphic to
$$
 \coker\bigl[
\bigl(j_{\bar D}^*\times \pi_D^*\bigr):  k^{\rm lc}(\bar X)^*\times k^{\rm lc}(D)^*\rightarrow 
k^{\rm lc}(\bar D)^*\bigr],
\eqno{(\ref{pwtlb.say}.2)}
$$
where $j_{\bar D} :\bar D\into \bar X$ is the natural injection.
Observe that  both $T_1:=k^{\rm lc}(\bar X)^*\times k^{\rm lc}(D)^*$ and $T_2:=k^{\rm lc}(\bar D)^*$
are $k$-tori, and then 
$$
\pic^{\rm pt}(X)\cong
\coker\bigl[T_1(k)\to T_2(k)\bigr].
\eqno{(\ref{pwtlb.say}.3)}
$$
(Note that although $\coker\bigl[T_1\to T_2\bigr]$ is a $k$-torus, usually we have a strict containment
$\coker\bigl[T_1(k)\to T_2(k)\bigr]\subsetneq \coker\bigl[T_1\to T_2\bigr](k)$,
so $\pic^{\rm pt}(X) $ is not the group of $k$-points of a $k$-torus.)

Let now $X$ be an arbitrary  proper, seminormal $k$-scheme.  
A {\it piecewise trivial divisorial sheaf} on $X$ is 
a divisorial sheaf on $X$ obtained as the push-forward of a 
piecewise trivial line bundle on its semi-regular locus $X^{\rm sr}$ (\ref{dn.notation}.3).  They form a group, for which we use either the divisor or the sheaf  notation
$$
\mcl^{\rm pt}(X):=\dsh^{\rm pt}(X):=\pic^{\rm pt}\bigl(X^{\rm sr}\bigr).
\eqno{(\ref{pwtlb.say}.4)}
$$

{\it Note.} Another variant of this definition would allow line bundles $L$
on $X^{\rm sr}$ that are assumed to be trivial only on $\bar X^{\rm sr}$. In this case  $L|_{D^{\rm sr}}$ need not be trivial, but it is a torsion element of $\pic(D^{\rm sr})$. This is similar to the usually minor  difference between  $\pic^{\circ}(X)$ and $\pic^{\tau}(X)$. 

\end{say}

\begin{say}[Proof of Theorem~\ref{main.thm}]
The surjective homomorphism
$$
\pi^*: \mcl(X) \onto \mcl(\bar X, \bar D)\cong  \oplus_i \cl(\bar X_i),
$$
was described in Paragraph~\ref{cwtlb.say}, 
proving (\ref{main.thm}.2). 

If $X$ is normal then $\mcl(X)\cong \cl(X)$ by Lemma~\ref{norm.mum=weil.lem} and the 
curve case is described in Paragraph~\ref{pwtlb.say.1}. Thus assume from now on that $\dim X\geq 2$ and $X$ is not normal.

Let $B=\{B_i:i\in I\}$ be a Mumford divisor such that $[B]\in \ker \pi^*$. 
 Then on each $\bar X_i$ we have a rational function $\bar\phi_i^B$ such that $(\bar\phi_i^B)=\pi_i^*B_i$. Together they give a rational function $\bar \phi^B$ on $\bar X$ that is regular and nonzero at all generic points of $\bar D$. 
Thus  $\bar \phi^B$ pulls back to a rational function  $\bar\phi^B_D$ on $\bar D^n$. 
If  $\bar\phi^B_D$ is the pull-back of a rational function  $\phi^B_D$ on $D^n$  then $\bar \phi^B$  descends to a
rational function $\phi^B$ on $X$ by (\ref{conductor.say}.1)  and  $(\phi^B)=B$.

Note that $\bar\phi^B$  and $\bar\phi^B_D$  are not uniquely determined by $B$, we can multiply them by constants. Thus the well-defined  object is the divisor
$(\bar\phi^B_D)$  of $\bar\phi^B_D$ on $\bar D^n$,  and
$$
(\bar\phi^B_D)\in  k\bigl( \bar D^n\bigr)^*/k^{\rm lc}\bigl( \bar D^n\bigr)^*\cong
\oplus_j k\bigl( \bar D^n_j\bigr)^*/k^{\rm lc}\bigl( \bar D^n_j\bigr)^*.
$$
Taking a further quotient by
$\pi_D^* k\bigl( D_j\bigr)^*$  on $\bar D_j$ we get  classes
$$
\bigl[\bar\phi^B_D\bigr]_j\in k\bigl( \bar D_j\bigr)^*/
\bigl\langle k\bigl( D_j\bigr)^*, k^{\rm lc}\bigl( \bar D_j\bigr)^*\bigr\rangle
$$
as  invariants of $[B]\in \ker \pi^*$. 
Since the restriction map $k(\bar X)\map k(\bar D)$ is surjective, 
these give  a  surjective map 
$$
\rho_D: \mcl^{\rm ct}(X)
\to \oplus_j \  k\bigl( \bar D^n_j\bigr)^*/
\bigl\langle k\bigl( D^n_j\bigr)^*, k^{\rm lc}\bigl( \bar D^n_j\bigr)^*\bigr\rangle.
$$
Once we prove that its kernel is $\mcl^{\rm pt}(X) $, we get the 
isomorphism in 
(\ref{main.thm}.5.a). (The definitions of  $k(D)$ and of  $k^{\rm lc}(D)$ are set up to be birational invariants of $D$.)

If  $\bar\phi^B_D$ is in the kernel of $\rho_D$ then
$\bar\phi^B_D|_{\bar D^n_j}=\psi_j\cdot \bar c_j$ where
$\psi_j\in k\bigl( D^n_j\bigr)^*$ and 
$\bar c_j\in  k^{\rm lc}\bigl( \bar D^n_j\bigr)^*$. Moreover, the product decomposition is unique up to a factor
$c_j\in  k^{\rm lc}\bigl( D^n_j\bigr)^*$.
Thus  $\psi_j$ determines a projective trivialization (\ref{trivialize.say}) of a line bundle
$L_j$ on $D^n_j$ whose pull-back to $\bar D^n_j$ is isomorphic to 
 $\o_{\bar X}(\bar B)|_{\bar D^n_j}$. Thus the kernel of $\rho_D$ consists of
piecewise trivial divisorial sheaves on $X$,  as described in
(\ref{pwtlb.say}).
 This shows (\ref{main.thm}.5.a) which  in turn implies (\ref{main.thm}.5.b--c) by the computations of Paragraph~\ref{mult.ff.say}.

Finally (\ref{main.thm}.5.b--c)  imply (\ref{main.thm}.4). 
\qed
\end{say}

\begin{say}[The multiplicative group of a function field] \label{mult.ff.say}
Let $V$ be a normal,  irreducible, proper
$k$-variety of dimension $\geq 1$ with field of constants $k^{\rm lc}(V)$. 
We saw in Definition~\ref{loc.const.say} that 
$$
k(V)^*/k^{\rm lc}(V)^*\cong \pdiv(V)\cong \z^{|\bar k|}.
\eqno{(\ref{mult.ff.say}.1)}
$$
Let next $W\to V$ be an irreducible,  separable  cover of proper, normal varieties.  For  any prime divisor $P\subset V$ let $P_1,\dots, P_{r(P)}$ be the prime divisors on $W$ lying over $P$. Write  $\pi^*P=\sum_i e_iP_i$ and set $e_P:=\gcd(e_1,\dots, e_{r(P)})$. Taking the divisor of a function gives
$$
q_{W/V}:k(W)^*/\langle k(V)^*, k^{\rm lc}(W)^*\rangle \to \sum_P \frac{\tsum_i \z[P_i]}{\tsum_i({e_i}/{e_P}) P_i}.
\eqno{(\ref{mult.ff.say}.2)}
$$
The kernel of $q_{W/V}$ consists of $\phi\in k(W)$ whose divisor $(\phi)$ is the pull-back of a $\q$-divisor on $V$. If $\deg(W/V)=m$  and $\phi\in \ker q$ then 
$$
\bigl(\phi^m\bigr)= \pi^*\pi_*(\phi)=\pi^*\bigl(\norm_{W/V}\phi\bigr).
\eqno{(\ref{mult.ff.say}.3)}
$$
In particular, 
the kernel of $q_{W/V}$  is $m$-torsion. 
Moreover, whether or not a $\q$-divisor on $V$ is a principal $\z$-divisor is invariant under separable field extensions. We can thus bound the torsion after base change to the separable closure $K\supset k$. In this case
we have a factorization
$W_K\to W^{\rm ab}_K\to V_K$ where $W^{\rm ab}_K\to V_K$ corresponds to the maximal abelian subfield of $K(W)\supset K(W^{\rm ab})\supset K(V)$. Then Kummer theory tells us that the
torsion in $K(W^*)/K(V)^*$ is  isomorphic to
$\gal\bigl(K(W^{\rm ab})/K(V)\bigr)$; cf. \cite[Sec.X.3]{MR554237}.

We can be especially explicit in the degree 2 case if $\chr k\neq 2$. 
Then  $k(W)/k(V)$  is a  Galois extension and there is a unique
$\psi\in k(W)^*/k(V)^*$ such that $\psi^2\in k(V)^*$. Thus, in this case
$$
k(W)^*/\langle k(V)^*, k^{\rm lc}(W)^*\rangle  \cong  
\left\{
\begin{array}{l}
\z^{|\bar k|}\qtq{if} k^{\rm lc}(W)\neq k^{\rm lc}(V)\qtq{and}\\
\z/2\z+ \z^{|\bar k|}\qtq{if} k^{\rm lc}(W)=k^{\rm lc}(V).
\end{array}
\right.
\eqno{(\ref{mult.ff.say}.4)}
$$
If $W\to V$ is purely inseparable of degree $q$ then
$k(W)^q\subset k(V)$, hence in this case
$$
k(W)^*/\langle k(V)^*, k[W]^*\rangle  \cong  (\mbox{$q$-torsion group of cardinality $|k|$}).
\eqno{(\ref{mult.ff.say}.5)}
$$
Finally an arbitrary $W\to V$ can be written as the composite of a
 purely inseparable map (say of degree $q$) and of a separable map (say of degree $m$)  and then  combining (\ref{mult.ff.say}.4--5) and elementary arguments as in \cite[p.880]{lang-alg} show that 
$$
k(W)^*/\langle k(V)^*, k[W]^*\rangle  \cong  (\mbox{finite $m$-torsion})+ (\mbox{$q$-torsion}) + \z^{|\bar k|}.
\eqno{(\ref{mult.ff.say}.6)}
$$
\end{say}

Riemann-Roch tells us that the  Euler characteristic of a   line bundle $L$ on a proper variety $X$ can be computed from it class $[L]\in \ns(X)$. 
The following example shows that for a Mumford divisor  $B$, 
its pull-back  $\bar B$ does not carry enough information to compute 
$\chi\bigl(X, \o_X(B)\bigr)$.

\begin{exmp}\label{dsh.on.dsm.say}
 Let $(X,D)$ be demi-smooth and $B$ a Mumford divisor on $X$.
Etale locally we can write $X$ as 
$$
X=(x_{11}x_{12}=0)\subset \a^{n+1}(x_{11}, x_{12}, x_2, \dots, x_n).
\eqno{(\ref{dsh.on.dsm.say}.1)}
$$
Its irreducible components are $X_1:=(x_{12}=0)$ and $X_2:=(x_{11}=0)$.
The divisors $B_i$ on $X_i$ are given by local equations
$B_i=\bigl( F_i( x_{1i}, x_2, \dots, x_n)=0\bigr)$ where the $F_i$ are rational functions. 
The sections of $\o_{X_i}(B_i)$ are locally of the form $G_i/F_i$ where $G_i$ is regular. 
Set 
$$
f_i:=F_i(0, x_2, \dots, x_n) \qtq{and} g_i:=G_i(0, x_2, \dots, x_n).
\eqno{(\ref{dsh.on.dsm.say}.2)}
$$
Two sections $ G_i/F_i$  patch to a section of $\o_X(B)$ iff
$$
g_1(x_2, \dots, x_n)/f_1(x_2, \dots, x_n)=
g_2(x_2, \dots, x_n)/f_2(x_2, \dots, x_n).
\eqno{(\ref{dsh.on.dsm.say}.3)}
$$
It is easier to express this in terms of divisors.
Write  $B_i|_D=\sum_j m_{ij}A_j$ where the $A_j$ are prime divisors on $D$.
Define
$$
\min\{B|_D\}:=\tsum_j \min\{m_{1j}, m_{2j}\}A_j.
\eqno{(\ref{dsh.on.dsm.say}.4)}
$$
Thus we get an exact sequence
$$
0\to \pi_*\o_{\bar X}(\bar B-\bar D)\to \o_X(B)\to
\o_D\bigl(\min \{B|_D\}\bigr) \to 0.
\eqno{(\ref{dsh.on.dsm.say}.5)}
$$

\medskip

{\it Claim \ref{dsh.on.dsm.say}.6.}  Fix  $X$ and $[\bar B]\in \cl(\bar X)$. Then
$\min \{B|_D\} $ can be arbitrary. 
\medskip

Proof. Pick any divisor $A$ on $D$. 
Choose $H$ sufficiently ample on $X$ such that
$H|_{X_i}-B_i$ and  $H|_D-A$ are very ample. 
Choose $H_i\sim H|_{X_i}$ such that  $H_i|_D=A+R_i$ and the $R_i$ are general. 
Then choose general  $G_i\sim H|_{X_i}-B_i$. Then
$B'_i:=H_i-G_i\sim B_i$ and 
$$
\begin{array}{rclcc}
\min \{B'|_D\}&=&\min \{(H_1-G_1)|_D, (H_2-G_2)|_D\}&&\\
&=&
\min \{A+R_1-G_1|_D, A+R_2-G_2|_D\}&=&A. \qed
\end{array}
$$

\medskip

Loosely speaking, this says that the group
$\mcl^{\rm ct}(X)/\mcl^{\rm pt}(X) $ is really important in cohomological questions.

\medskip
{\it Informal corollary \ref{dsh.on.dsm.say}.7.} 
Let $X$ be a proper, demi-normal, non-normal variety  and $B$ a Mumford divisor on $X$.
Knowing  $X$ and $[\bar B]\in \cl(\bar X)$ tells essentially nothing about the Euler characteristic or the cohomology groups of   $\o_X(B)$. \qed

\medskip
{\it Informal corollary \ref{dsh.on.dsm.say}.8.} Let $X$ be a proper, demi-normal, non-normal variety of dimension $n$,  $B$ a Mumford divisor on $X$ and $H$ an ample Cartier divisor. Then
 $X$ and $[\bar B]\in \cl(\bar X)$ determine  the 2 highest coefficients 
of the Hilbert polynomial
$$
\chi\bigl(X, \o_X(B+mH)\bigr),
$$
but give no information about the others.  \qed

\end{exmp} 

\begin{ack} I thank   A.J.\ de~Jong and Chenyang~Xu for helpful  comments. 
Partial  financial support    was provided  by  the NSF under grant number
 DMS-1362960.
\end{ack}


\begin{thebibliography}{BVRS09}

\bibitem[AB69]{MR0266923}
A.~Andreotti and E.~Bombieri, \emph{Sugli omeomorfismi delle variet\`a
  algebriche}, Ann. Scuola Norm. Sup Pisa (3) \textbf{23} (1969), 431--450.
  \MR{0266923 (42 \#1825)}

\bibitem[ACG11]{acg2}
Enrico Arbarello, Maurizio Cornalba, and Pillip~A. Griffiths, \emph{Geometry of
  algebraic curves. {V}olume {II}}, Grundlehren der Mathematischen
  Wissenschaften, vol. 268, Springer, Heidelberg, 2011, With a contribution by
  Joseph Daniel Harris. \MR{2807457 (2012e:14059)}

\bibitem[Alb32]{alb-20}
Giacomo Albanese, \emph{Corrispondenze algebriche fra i punti di due superfizie
  algebriche}, Boll. Uni. Math. Ital. \textbf{11} (1932), 131--138.

\bibitem[Alb34]{alb-21}
\bysame, \emph{Corrispondenze algebriche fra i punti di due superfizie
  algebriche}, Ann. Sci. Norm. Sup. Pisa. \textbf{3} (1934), 1--26, 149--182.

\bibitem[Alb96]{MR1419986}
\bysame, \emph{Collected papers of {G}iacomo {A}lbanese}, Queen's Papers in
  Pure and Applied Mathematics, vol. 103, Queen's University, Kingston, ON,
  1996, With a biography of Albanese in English by Ciro Ciliberto and Edoardo
  Sernesi, Edited by Ciliberto, Paulo Ribenboim and Sernesi. \MR{1419986}

\bibitem[Amb03]{ambro}
Florin Ambro, \emph{Quasi-log varieties}, Tr. Mat. Inst. Steklova \textbf{240}
  (2003), no.~Biratsion. Geom. Linein. Sist. Konechno Porozhdennye Algebry,
  220--239. \MR{1993751 (2004f:14027)}

\bibitem[AN67]{MR0239118}
Aldo Andreotti and Fran{\cedilla{c}}ois Norguet, \emph{La convexit\'e
  holomorphe dans l'espace analytique des cycles d'une vari\'et\'e
  alg\'ebrique}, Ann. Scuola Norm. Sup. Pisa (3) \textbf{21} (1967), 31--82.
  \MR{0239118 (39 \#477)}

\bibitem[BN12]{MR3044492}
John Brevik and Scott Nollet, \emph{Picard groups of normal surfaces}, J.
  Singul. \textbf{4} (2012), 154--170. \MR{3044492}

\bibitem[BVRS09]{MR2457299}
Luca Barbieri-Viale, Andreas Rosenschon, and V.~Srinivas, \emph{The
  {N}\'eron-{S}everi group of a proper seminormal complex variety}, Math. Z.
  \textbf{261} (2009), no.~2, 261--276. \MR{2457299}

\bibitem[BVS93]{MR1203911}
L.~Barbieri~Viale and V.~Srinivas, \emph{On the {N}\'eron-{S}everi group of a
  singular variety}, J. Reine Angew. Math. \textbf{435} (1993), 65--82.
  \MR{1203911}

\bibitem[dFKX12]{dkx}
Tommaso de~Fernex, J{\'a}nos Koll{\'a}r, and Chenyang Xu, \emph{{The dual
  complex of singularities}}, ArXiv e-prints (2012).

\bibitem[Fuj17]{fuj-book}
Osamu Fujino, \emph{Foundations of the minimal model program}, MSJ Memoirs,
  vol.~35, Mathematical Society of Japan, Tokyo, 2017. \MR{3643725}

\bibitem[Ful84]{Fulton84}
William Fulton, \emph{Intersection theory}, Ergebnisse der Mathematik und ihrer
  Grenzgebiete (3), vol.~2, Springer-Verlag, Berlin, 1984. \MR{MR732620
  (85k:14004)}

\bibitem[Gro62]{FGA}
Alexander Grothendieck, \emph{Fondements de la g\'eom\'etrie alg\'ebrique.
  [{E}xtraits du {S}\'eminaire {B}ourbaki, 1957--1962.]}, Secr\'etariat
  math\'ematique, Paris, 1962. \MR{MR0146040 (26 \#3566)}

\bibitem[Har94]{MR1291023}
Robin Hartshorne, \emph{Generalized divisors on {G}orenstein schemes},
  Proceedings of Conference on Algebraic Geometry and Ring Theory in honor of
  Michael Artin, Part III (Antwerp, 1992), vol.~8, 1994, pp.~287--339.
  \MR{MR1291023 (95k:14008)}

\bibitem[HP15]{MR3356940}
Robin Hartshorne and Claudia Polini, \emph{Divisor class groups of singular
  surfaces}, Trans. Amer. Math. Soc. \textbf{367} (2015), no.~9, 6357--6385.
  \MR{3356940}

\bibitem[Jaf88]{MR975125}
David~B. Jaffe, \emph{Space curves which are the intersection of a cone with
  another surface}, Duke Math. J. \textbf{57} (1988), no.~3, 859--876.
  \MR{975125}

\bibitem[KM98]{km-book}
J{\'a}nos Koll{\'a}r and Shigefumi Mori, \emph{Birational geometry of algebraic
  varieties}, Cambridge Tracts in Mathematics, vol. 134, Cambridge University
  Press, Cambridge, 1998, With the collaboration of C. H. Clemens and A. Corti,
  Translated from the 1998 Japanese original.

\bibitem[Kol96]{rc-book}
J{\'a}nos Koll{\'a}r, \emph{Rational curves on algebraic varieties}, Ergebnisse
  der Mathematik und ihrer Grenzgebiete. 3. Folge., vol.~32, Springer-Verlag,
  Berlin, 1996.

\bibitem[Kol13]{kk-singbook}
\bysame, \emph{Singularities of the minimal model program}, Cambridge Tracts in
  Mathematics, vol. 200, Cambridge University Press, Cambridge, 2013, With the
  collaboration of S{\'a}ndor Kov{\'a}cs.

\bibitem[Kol17]{k-coherent}
\bysame, \emph{Coherence of local and global hulls}, Methods Appl. Anal.
  \textbf{24} (2017), no.~1, 63--70. \MR{3694300}

\bibitem[Kol18a]{k-lpg1}
\bysame, \emph{{Log-plurigenera in stable families}}, ArXiv e-prints (2018).

\bibitem[Kol18b]{k-lpg2}
\bysame, \emph{{Log-plurigenera in stable families of surfaces}}, ArXiv
  e-prints (2018).

\bibitem[Lan02]{lang-alg}
Serge Lang, \emph{Algebra}, third ed., Graduate Texts in Mathematics, vol. 211,
  Springer-Verlag, New York, 2002. \MR{1878556}

\bibitem[Mat52]{mats-pic}
Teruhisa Matsusaka, \emph{On the algebraic construction of the {P}icard
  variety. {II}}, Jap. J. Math. \textbf{22} (1952), 51--62. \MR{0062471}

\bibitem[N{\'e}r52]{neron-pic}
Andr\'e N{\'e}ron, \emph{La th\'eorie de la base pour les diviseurs sur les
  vari\'et\'es alg\'ebriques}, Deuxi\`eme {C}olloque de {G}\'eom\'etrie
  {A}lg\'ebrique, {L}i\`ege, 1952, Georges Thone, Li\`ege; Masson \& Cie,
  Paris, 1952, pp.~119--126. \MR{0052154}

\bibitem[Pic95]{picard-95}
\'Emile Picard, \emph{Sur la th\'eorie des groupes et des surfaces
  alg\'ebriques}, Rend. Circ. Math. Palermo \textbf{9} (1895).

\bibitem[Ros54]{MR0061422}
Maxwell Rosenlicht, \emph{Generalized {J}acobian varieties}, Ann. of Math. (2)
  \textbf{59} (1954), 505--530. \MR{0061422}

\bibitem[RS09]{MR2567426}
G.~V. Ravindra and V.~Srinivas, \emph{The {N}oether-{L}efschetz theorem for the
  divisor class group}, J. Algebra \textbf{322} (2009), no.~9, 3373--3391.
  \MR{2567426}

\bibitem[Ser59]{MR0103191}
Jean-Pierre Serre, \emph{Groupes alg\'ebriques et corps de classes},
  Publications de l'institut de math\'ematique de l'universit\'e de Nancago,
  VII. Hermann, Paris, 1959. \MR{0103191}

\bibitem[Ser79]{MR554237}
\bysame, \emph{Local fields}, Graduate Texts in Mathematics, vol.~67,
  Springer-Verlag, New York-Berlin, 1979, Translated from the French by Marvin
  Jay Greenberg. \MR{554237}

\bibitem[Sev1906]{MR1511373}
Francesco Severi, \emph{Sulla totalit\`a delle curve algebriche tracciate sopra
  una superficie algebrica}, Math. Ann. \textbf{62} (1906), no.~2, 194--225.
  \MR{1511373}

\bibitem[Sev1913]{severi-13}
\bysame, \emph{Una teorema di inversione per gli imtegrali semplici $1^{\rm a}$
  specie appartementi ad una superficie algebriche}, Atti. Inst. Veneto
  \textbf{72} (1913), 765--772.

\bibitem[Sev47]{MR0024985}
\bysame, \emph{Funzioni quasi abeliane}, Pontificiae Academiae Scientiarum
  Scripta Varia, v. 4, Publisher unknown, Vatican City, 1947. \MR{0024985}

\bibitem[Sho92]{sho-3ff}
V.~V. Shokurov, \emph{Three-dimensional log perestroikas}, Izv. Ross. Akad.
  Nauk Ser. Mat. \textbf{56} (1992), no.~1, 105--203. \MR{1162635 (93j:14012)}

\bibitem[{Sta}15]{stacks-project}
The {Stacks Project Authors}, \emph{{S}tacks {P}roject},
  http://stacks.math.columbia.edu, 2015.

\bibitem[Tra70]{MR0277542}
Carlo Traverso, \emph{Seminormality and {P}icard group}, Ann. Scuola Norm. Sup.
  Pisa (3) \textbf{24} (1970), 585--595. \MR{0277542 (43 \#3275)}

\bibitem[Zar71]{MR0469915}
Oscar Zariski, \emph{Algebraic surfaces}, supplemented ed., Springer-Verlag,
  New York-Heidelberg, 1971, With appendices by S. S. Abhyankar, J. Lipman, and
  D. Mumford, Ergebnisse der Mathematik und ihrer Grenzgebiete, Band 61.
  \MR{0469915}

\end{thebibliography}

\def\cprime{$'$} \def\cprime{$'$} \def\cprime{$'$} \def\cprime{$'$}
  \def\cprime{$'$} \def\cprime{$'$} \def\dbar{\leavevmode\hbox to
  0pt{\hskip.2ex \accent"16\hss}d} \def\cprime{$'$} \def\cprime{$'$}
  \def\polhk#1{\setbox0=\hbox{#1}{\ooalign{\hidewidth
  \lower1.5ex\hbox{`}\hidewidth\crcr\unhbox0}}} \def\cprime{$'$}
  \def\cprime{$'$} \def\cprime{$'$} \def\cprime{$'$}
  \def\polhk#1{\setbox0=\hbox{#1}{\ooalign{\hidewidth
  \lower1.5ex\hbox{`}\hidewidth\crcr\unhbox0}}} \def\cdprime{$''$}
  \def\cprime{$'$} \def\cprime{$'$} \def\cprime{$'$} \def\cprime{$'$}
\providecommand{\bysame}{\leavevmode\hbox to3em{\hrulefill}\thinspace}
\providecommand{\MR}{\relax\ifhmode\unskip\space\fi MR }
\providecommand{\MRhref}[2]{%
  \href{http://www.ams.org/mathscinet-getitem?mr=#1}{#2}
}
\providecommand{\href}[2]{#2}

\bigskip

\noindent  Princeton University, Princeton NJ 08544-1000

{\begin{verbatim} kollar@math.princeton.edu\end{verbatim}}

\end{document}